\documentclass[12pt]{article}
 \usepackage{latexsym}
 \usepackage{graphicx}
 \usepackage{amssymb}
 \title{Edge Decompositions of Hypercubes by Paths and by Cycles}

 \author{Michel Mollard\thanks{CNRS Universit\'e Joseph Fourier}\\
 \small Institut Fourier \\[-0.8ex]
 \small 100, rue des Maths\\[-0.8ex]
 \small 38402 St Martin d'H\`eres Cedex FRANCE \\
 \small \texttt{michel.mollard@ujf-grenoble.fr}\\
 \small and\\
 Mark Ramras\\
 \small Department of Mathematics\\[-0.8ex]
 \small Northeastern University\\[-0.8ex]
  \small Boston, MA 02115, USA\\[-0.8ex]
 \small \texttt{m.ramras@neu.edu}}
 \begin{document}
 \maketitle
 \date{}
 \begin{abstract}
 If $H$ is isomorphic to a subgraph of $G$, we say that $H$
 {\it divides} $G$ if there exist embeddings $\theta_1, \theta_2, \ldots , \theta_k$ of $H$ such that  
 $$\{\{E(\theta_1(H)), E(\theta_2(H)), \ldots, E(\theta_k(H))\}$$
  is a partition of $E(G)$.  For purposes of simplification we will often omit the embeddings, saying that we have an edge decomposition by copies of $E(H)$.  
 
 Many authors  
  have studied this notion for various subgraphs of
 hypercubes.  We continue such a study in this paper.
 \end{abstract}
 
 \section{Introduction and Preliminary Results}
 \newtheorem{definition}{Definition}
 \begin{definition}
 If $H$ is isomorphic to a subgraph of $G$, we say that $H$
 {\it divides} $G$ if there exist embeddings $\theta_1, \theta_2, \ldots , \theta_k$ of $H$ such that 
 $$\{\{E(\theta_1(H)), E(\theta_2(H)), \ldots, E(\theta_k(H))\}$$
  is a partition of $E(G)$. 
 \end{definition}
 
 Ramras \cite{Ram} has defined a more restrictive concept.
 \begin{definition}
 A {\it fundamental set} of edges of a graph $G$ is a subset of $E(G)$
 whose translates under some subgroup of the automorphism group of $G$
 partition $E(G)$.
 \end{definition}
 
 Edge decompositions of graphs by subgraphs have a long history.  For example, there is a Steiner triple system of order $n$ if and only if the complete graph $K_n$ has an edge-decomposition 
 by $K_3$.  In 1847 Kirkman \cite{Kirk1} proved that for a Steiner triple system to exist it is necessary that $n\equiv 1 \pmod {6}$ or $n\equiv 3 \pmod {6}$.  In 1850 he  proved the converse holds also \cite{Kirk2}.
 \newtheorem{theorem}{Theorem}\label{Kirk1} \label{Kirk2}
 \begin{theorem}  A Steiner system of order $n\geq 3$ exists if and only if $n\equiv 1 \pmod {6})$ or $n\equiv 3 \pmod {6}$.
 \end{theorem}
 
 In more modern times (1964) G. Ringel \cite{Ring2} stated the following conjecture, which is still open.
 \newtheorem{conjecture}{Conjecture}
  \begin{conjecture}
 If $T$ is a fixed tree with $m$ edges then $K_{2m+1}$ is edge-decomposable into $2m+1$ copies of $T$.
 \end{conjecture}
 
 By $Q_n$ we mean the $n$-dimensional hypercube.  We regard its vertex set, $V(Q_n)$, as $\mathcal P(\{1, 2, \ldots, n\})$, the set of subsets of $\{1, 2, \ldots , n\}$.
 Two vertices $x$ and $y$ are considered adjacent (so $\langle x, y\rangle \in E(Q_n)$) if $|\,x\Delta y\,|=1$, where $\Delta$ denotes the symmetric difference of the two subsets $x$ and $y$.  $(V(Q_n), \Delta)$ is isomorphic as a group to $(\mathbb{Z}_2^n, +)$.  Occasionally, when convenient, we shall use the vector notation for vertices; thus $\vec{x}$ and $\vec{y}$ are adjacent precisely when they differ in exactly one component.
 Note that for $k<n, \mathcal P(\{1, 2, \dots, k\}) \subset \mathcal P(\{1, 2, \dots, n\})$ so that $V(Q_k)\subset V(Q_n)$.  In fact, from the definition of adjacency, it follows that $Q_k$ is an induced subgraph
 of $Q_n$.
 
 Beginning in the early 1980's, interest in hypercubes (and similar hypercube-like networks such as ``cube-connected cycles" and ``butterfly" networks) increased dramatically with the construction of massively parallel-processing computers, such as the ``Connection Machine" whose architecture is that of the $16$-dimensional hypercube, with $2^{16}=65,536$ processors as the vertices.  Problems of routing message packets simultaneously along paths from one processor to another led to an interest in questions of edge decompositions of $E(Q_n)$
 by paths.  An encyclopedic discussion of this and much more can be found in \cite{Lei}.
 
 In \cite{Ram} we have shown that if ${\cal G}$ is a subgroup of ${\rm
Aut}(Q_n)$
 and for all $g\in {\cal G}$, with $g\neq  id$ (where $id$ denotes the identity
 element), $g(E(H))\cap E(H)=\emptyset$, then there is a packing of these
 translates of $E(H)$ in $Q_n$, {\it i.e.} they are pairwise disjoint.  If, in
 addition, $|\,E(H)\,|\cdot |\,{\cal G}\,|=n\cdot 2^{n-1}=|\,E(Q_n)\,|$, then
the
 translates of $E(G)$ by
 the elements of ${\cal G}$ yield an edge decomposition of $Q_n$.  In
 \cite{Ram} it is shown that every tree on $n$ edges can be embedded in
 $Q_n$ as a fundamental set.  (This result for edge decompositions was obtained independently by Fink \cite{Fink}). In \cite{Ram2} this is extended to
 certain trees
 and certain cycles on $2n$ edges.  Decompositions of $Q_n$ by
 $k$-stars are proved for all $k\leq n$ in \cite{etal}.  Recently,
 Wagner and Wild \cite{WagWil} have constructed, for each value of $n$, a tree
on $2^{n-1}$
 edges that is a fundamental
 set for $Q_n$.
The structure of ${\rm Aut}(Q_n)$ is discussed in \cite{Ram}.  For
 each subset $A$ of $\{1,2,\ldots ,n\}$, the complementing
 automorphism $\sigma_A$ is defined by
 $\sigma_A(x)=A\Delta \{x\}.$   
 Another type of automorphism arises from the group of permutations
 ${\cal S}_n$ of $\{1,2,\ldots, n\}$.  For $x=\{x_1, x_2, \ldots, x_m\} \subseteq \{1, 2, \ldots, n\}$ and $\theta \in {\cal S}_n$ we denote by $\rho_{\theta(x)}$ the vertex $\{\theta(x_1), \theta(x_2), \ldots , \theta(x_m)\}$.
  The mapping $\rho_{\theta}: V(Q_n)\longrightarrow V(Q_n)$ defined in this way is easily seen to belong to ${\rm Aut}(Q_n)$.  Every
 automorphism in ${\rm Aut}(Q_n)$ can be expressed uniquely in the form
 $\sigma_A\circ \rho_{\theta}$, where this notation means that we first apply $\rho_{\theta}$.  Note:  $\rho_{\theta}\circ \sigma_A=\sigma_{\theta(A)}\circ \rho_{\theta}$.
 
 To avoid ambiguity in what follows we make this definition:
  \begin{definition}

  By $P_k$, the  ``$k$-path'', we mean the path with $k$ edges.
  \end{definition}
 
 \noindent {\bf Questions}\\
 (1)  For which $k$ dividing $n\cdot 2^{n-1}$ does $P_k$
 divide $Q_n$?\\
 (2)  For which $k$ dividing $n\cdot 2^{n-1}$ does $C_k$, the cycle on $k$
edges,
 divide $Q_n$?\\
 (3)  For those $k$ for which the answer to either (1) or (2) is
 ``yes", is the edge set used in the decomposition a fundamental set for $Q_n$?
 
 We begin this introductory section with some examples.  In later
 sections we prove a variety of results relating to these questions,
 and in the final section we summarize our findings.
 
 \newtheorem{example}{Example}
 \begin{example}
 \end{example}
 
 Let $T$ be the 2-star ($=$ the $2$-path) contained in $Q_3$ with
 center $000$, and leaves $100,
  010$.  Then ${\cal G}=\{id, \sigma_{123}, \sigma_1\rho_{(123)},
 \sigma_{12}\rho_{(132)},
 \sigma_3\rho_{(132)}, \sigma_{23}\rho_{(123)}\}$ is a (cyclic) subgroup of
 Aut$(Q_3)$ of order 6, and the 6 translates of $T$ under ${\cal G}$
 yield an edge
 decomposition of $Q_3$.  \hfill  $\Box$
 
 \vspace{.2in}
 
 Note, however, that ${\cal G}$ does not work for the 2-star $T^{\prime}$, whose
 center is $000$ and whose leaves are $100$ and $001$.  The subgroup which
 works for this 2-star is ${\cal G}^{\prime}=\{id, \sigma_{123},
 \sigma_1\rho_{(132)},
 \sigma_{13}\rho_{(123)}, \sigma_2\rho_{(123)}, \sigma_{23}\rho_{(132)}\}$.\\
 
 \begin{example}
 \end{example}
 
 $P_6$ does not divide $Q_3$.  For since $Q_3$ has 12 edges, if $P_6$ {\it did}
 divide $Q_3$ then $Q_3$ would have an edge-decomposition consisting of 2 copies
 of $P_6$.  The degree sequence (in decreasing order) of each $P_6$ is $2,2,2,2,2,1,1,0$, whereas $Q_3$, of
 course, is 3-regular.  Thus the vertex of degree $0$ in one $P_6$ would require a degree of $3$ in the other, which is impossible.  \hfill  $\Box$

 \begin{example}
 \end{example}
 
 $P_4$ does not divide $Q_3$.  Since $P_4$ has 4 edges, we would need 3 copies
 of $P_4$ for an edge-decomposition of $Q_3$. Call the three copies of $P_4$ $P^{(1)}, P^{(2)},$ and
 $P^{(3)}$.  At each vertex $v$ of $Q_3$, $\sum_{1\leq i\leq 3}
 {\rm deg}_{P^{(i)}}(v)
 =3$. Label the vertices of $Q_3$ $(v_1)$ to $(v_8)$ such that the degree sequence of $P^{(1)}$, is decreasing. Consider the $3\times 8$ array ${\rm deg}_{P^(i)}(v_j)$ . The first row  is thus 2 2 2 1 1 0 0 0.
 In the second and third rows, in order for the column sums to be 3, there must be exactly 3 1's (and 3 0's) in the first 3 columns.  Similarly, in
 the last 3 columns there must be exactly 3 1's (and 3 0's).  Thus in the second and third rows we have at least 6 1's, and so at least one of these
 rows must have at least 3 1's.  But each row is a permutation of the first, which has only 2 1's.  Contradiction.  Hence $P_4$ does not divide $Q_3$.  \hfill 
$\Box$

 \begin{example}
 \end{example}
 
 Since $Q_3$ is $3$-regular, the $4$-star is not a subgraph.
 The {\it other} tree on 4 edges {\em does} divide $Q_3$.  Let $T$ be the 3-star
 centered at $000$ union the edge $\langle001, 101\rangle$.  Let ${\cal
 G}=<\sigma_{23}\rho_{(123)}>$,
 which is a cyclic subgroup of Aut$(Q_3)$ of order 3.  A
 straight-forward calculation shows that the translates of $T$ under
 ${\cal G}$ form an edge decomposition of $Q_3$.  \hfill  $\Box$

 \newtheorem{proposition}{Proposition}
 \begin{proposition}  \label{doesNotdiv}
 For $k\geq 3, P_{2^k}$ does not divide $Q_{2k+1}$.
 \end{proposition}
 
  \noindent {\em Proof.}     Suppose that $k\geq 3$, and suppose that $P_{2^k}$
 divides $Q_{2k+1}$. The matrix $\left(a_{iv}\right)$ formed by the degree sequences of copies of $P_{2^k}$has $2^{2k+1}$ columns,
 and
 $$(2k+1)\cdot 2^{2k}/2^k=(2k+1)2^k$$
  rows. Then since each row has exactly two 1's, the entire matrix has
 \linebreak $(2k+1)2^{k+1}$ 1's.  But since each vertex of $Q_{2k+1}$
 has degree $2k+1$, each column sum is $2k+1$, and thus each column has
 at least one 1.  Thus there must be at least $2^{2k+1}$ 1's in the
 matrix.  Therefore, $(2k+1)2^{k+1}\geq 2^{2k+1}$.  This is equivalent
 to $2k+1\geq 2^{k}$.  But for $k\geq 3$ this is clearly false.  Thus
 for $k\geq 3, P_{2^k}$ does not divide $Q_{2k+1}$.
  \hfill  $\Box$

We will prove in Section 3 that for  $k=2, P_{2^k}$ \emph{does} divide $Q_{2k+1}$.

\vspace{.2in}
The next result is Proposition 8 of \cite{Ram2}.

 \begin{proposition}\label{Ram2, Prop8}
 Let $n$ be odd, and suppose that $P_k$ divides $Q_n$.  Then $k\leq n$.
 \end{proposition}
 
  \newtheorem{lemma}{Lemma}
 \begin{lemma} \label{Lemmadiv}
 ``Divisibility" is transitive, {\em i.e.} if $G_1$ divides $G_2$ and
 $G_2$ divides $G_3$, then $G_1$ divides $G_3$.
 \end{lemma}
 
 \noindent {\em Proof.}  This follow immediately from the definition of
 ``divides".  \hfill  $\Box$
 
 \newtheorem{corollary}{Corollary}
 \begin{corollary}\label{PkdivQn}
 If $k$ divides $n$ then $P_k$ divides $Q_n$.
 \end{corollary}
 
 \noindent {\em Proof.}  By \cite{Ram}, Theorem 2.3, $T$ divides
 $Q_n$ for {\em every} tree $T$ on $n$ edges.  In particular,
 then, $P_n$ divides $Q_n$.  Clearly, if $k$ divides $n$ then $P_k$
 divides $P_n$.  Hence, by Lemma \ref{Lemmadiv}, $P_k$ divides $Q_n$.
 \hfill $\Box$
 
 \vspace{.2in}
 
 We have the following partial converse.

 \begin{proposition} \label{kdivn}
 If $P_k$ divides $Q_n$ and $k$ is odd, then $k$ divides $n$.
 \end{proposition}
 
 \noindent {\em Proof.}  Since $P_k$ divides $Q_n$, $k$ divides $n\cdot
 2^{n-1}$.  But since $k$ is odd, this means that $k$ divides $n$.
 \hfill $\Box$
 
 \begin{definition}
 If $G_1$ and $G_2$ are graphs then by $G_1\Box G_2$ we mean the graph that is the Cartesian product of $G_1$ and $G_2$.
 \end{definition}
 
 \begin{lemma}\label{Hdiv G1 and G2}
 If $H$ divides $G_1$ and $H$ divides $G_2$ then $H$ divides $G_1\Box G_2$.
 \end{lemma}
 
 \noindent {\em Proof.}  This is obvious because $E(G_1\Box G_2)$ consists of $|\,V(G_1)\,|$ copies of $E(G_2)$ and $|\,V(G_2)\,|$ copies of $E(G_1)$.  \hfill  $\Box$

 \vspace{.2in}
 
 \begin{proposition}\label{QkdivQn}
 If $k$ divides $n$ then $Q_k$ divides $Q_n$.
 \end{proposition}
 
 \noindent {\em Proof.}  Let $n=mk$.  We argue by induction on $m$.  The statement is obvious for $m=1$.  Now let $m>1$ and assume the statement is true for $m-1$.  The desired result follows from
 Lemma \ref{Hdiv G1 and G2} and the fact that $Q_{(m-1)k}\Box Q_k\simeq Q_{(m-1)k+k}=Q_{mk}$.  \hfill  $\Box$

 The converse to Proposition \ref{QkdivQn} follows easily from the next lemma.
 
  \begin{lemma}\label{reg}
 Suppose that the subgraph $H$ of $G$ edge-divides $G$.  If $G$ is $n$-regular and $H$ is $k$-regular, then $k$ divides $n$.
 \end{lemma}
 
 \noindent {\em Proof.}  Since the copies of $E(H)$ form an edge-partition of $E(G)$, each vertex $v$ of $H$ must belong to exactly $n/k$ copies of $H$ and so $k$ divides $n$.  \hfill   $\Box$ 
 
 \begin{corollary}\label{Co: conversekdivn}
 If $Q_k$ divides $Q_n$ then $k$ divides $n$. 
 \end{corollary}
 
 \noindent {\em Proof.}  Since $Q_k$ is $k$-regular and $Q_n$ is $n$-regular, this follows immediately from Lemma \ref{reg}.  \hfill   $\Box$
 \vspace{.1in}
 
 Combining Proposition \ref{QkdivQn} and Corollary \ref{Co: conversekdivn} we obtain
 
 \begin{proposition}\label{QkdivQniff kdivn}
 $Q_k$ divides $Q_n$ if and only if $k$ divides $n$.
 \end{proposition}
 
 As an immediate consequence of Lemma \ref{Lemmadiv} and Proposition
 \ref{QkdivQn} we have
 
 \begin{corollary} \label{ImmedConseq}
 If $k$ divides $n$ and if $P_j$ divides $Q_k$ then $P_j$ divides $Q_n$.
 \end{corollary}
 
 We have a more general consequence.
 
 \begin{corollary}\label{TdivQn}
  If $k$ divides $n$ and $T$ is any tree on $k$ edges, then there is an embedding of $T$ which divides
 $Q_n$.  
 \end{corollary}
 
 \noindent {\em Proof.}  By \cite{Ram}, Theorem 2.3, by mapping any given vertex of $T$ to $\emptyset$ and assigning distinct labels $1, 2, \ldots, k$ to the edges of $T$ we get a subtree of $Q_k$ isomorphic to $T$ that divides $Q_k$.   Hence
 by Lemma \ref{Lemmadiv} and Proposition \ref{QkdivQn}, $T$ divides $Q_n$.  \hfill  $\Box$

 \begin{proposition}\label{Prop:  P_2^j}
 If $n$ is even, and $j<n$ then $P_{2^j}$ divides $Q_n$.
 \end{proposition}
 
 \noindent {\em Proof.}  It is proved in \cite{AS} that the cycle
 $C_{2^n}$ divides $Q_n$.  The Hamiltonian cycle $C_{2^n}$ is divisible
 by any path $P_q$, as long as $q$ divides $2^n$ and $q<2^n$.
 Thus $C_{2^n}$ is divisible by $P_{2^j}$ provided $j<n$.  The result
 now follows from Lemma \ref{Lemmadiv}.  \hfill  $\Box$
 
 \begin{proposition}\label{cycle}
 If $n$ is even, and $C$ is the $2n$-cycle with initial vertex $\emptyset$, and
 edge direction sequence $(1,2,\ldots , n)^2\stackrel{\rm def}\equiv(1,2,\ldots , n, 1, 2, \ldots, n)$, then $Q_n$ is edge-decomposed by
 the copies of $C$ under the action of
 ${\cal G}=\{\sigma_A\,|\,A\subset \{1, 2, \ldots, n-1\}$, $|\,A\,|$ even$\}$.  So $E(C)$ is
 fundamental for $Q_n$.
 \end{proposition}
 
 \noindent {\em Proof.}  $C$ consists of the path $P$, followed by
 $\sigma_{\{1,2,\ldots, n\}}(P)$, where $P$ is the path with initial
 vertex $\emptyset$ and edge direction sequence
 $1,2,\ldots , n$.  Note that for any $B\subseteq \{1, 2, \ldots,
 n\},$ for any edge $e, \sigma_B(e)=e$ implies that $ B=\emptyset$ or
 $|\,B\,|=1.$ Now we shall show that for every subset $A\subset \{1,2,\ldots,
 n-1\}$ with
 $|A|$ even, $\sigma_A(C)\cap C=\emptyset$.  It should be noted that
 these $A$'s form a subgroup
 of ${\rm Aut}(Q_n)$ of order $2^{n-2}$.  So suppose that $e=\langle
 x,y\rangle \in C\cap \sigma_A(C)$.  Let the direction of $e$ be $i$.
 Then the direction of $\sigma_A(e)$ is $i$.  If $A\neq \emptyset$, then
 since
 $|A|$ is even, $\sigma_A(e)\neq e$.  The only other edge in $C$ with
 direction $i$ is $\sigma_{\{1,2,\ldots , n\}}(e)$.  So if
 $\sigma_A(e)\in C$, then
 $\sigma_A(e)=\sigma_{\{1,2,\ldots ,n \}}(e)$.  Therefore $\sigma_A\cdot
 \sigma_{\{1,2,\ldots , n\}}(e)=e$, {\it i.e.} $\sigma_{A\Delta \{1,2,\ldots ,
 n\}}(e)=e$.  Since $A$ and $\{1,2,\ldots , n\}$
  are even, so is $A\Delta \{1,2,\ldots , n\}=\overline{A}$.  Hence
 $A\Delta \{1,2,\ldots , n\}=\emptyset$, {\it i.e.} $A=\{1,2,\ldots ,
 n\}$.  But $n\notin A$, so we have a contradiction.
 
 Thus we have a group ${\cal G}$ of automorphisms of $C$ of order
 $2^{n-2}$, such that for $g\in {\cal G}, g\neq id$ , $g(E(C))\cap
 E(C)=\emptyset$.  Furthermore, since $|E(C)|=2n$, it follows that
 $|{\cal G}|\cdot |E(C)|=|E(Q_n)|$.  Hence by \cite{Ram}, Lemma 1.1, the
 translates of $E(C)$ via the elements of ${\cal G}$ form an edge
 decomposition of $Q_n$.  \hfill  $\Box$

  \begin{corollary}\label{P2kdivQn}
 If $n$ is even, $k<n$ and $k$ divides $n$, then $P_{2k}$ divides $Q_n$.
 \end{corollary}
 
 \noindent {\em Proof.}  Since $k$ divides $n$, $2k$ divides $2n$, and
 thus since $2k<2n, P_{2k}$ divides the $2n$-cycle $C$ of  Proposition
 \ref{cycle}.  Hence by Proposition \ref{cycle}, $P_{2k}$ divides
 $Q_n$.  \hfill  $\Box$
 
 \begin{corollary}
 If $n$ and $k$ are both even and $k$ divides $n$, and $C$ is the $2k$-cycle with
 initial vertex $\emptyset$, and edge direction sequence $(1,2,\ldots ,
 k)^2$, then $C$ divides $Q_n$.
 \end{corollary}
 
 \noindent {\em Proof.}  By the proposition, $C$ divides $Q_k$, and by
 Proposition \ref{QkdivQn}, $Q_k$ divides $Q_n$.  The result now follows from
 Lemma \ref{Lemmadiv}.  \hfill  $\Box$

 \section{$P_4$ divides $Q_5$}
 
If $k$ is odd then by Proposition \ref{kdivn} and Lemma \ref{PkdivQn} $P_k$
divides $Q_n$ and only if $k$ divides $n$.  Thus the smallest value of $k$ for
which Question (1) remains open is $k=4.$  Corollary \ref{P2kdivQn} settles the
matter in the affirmative when $n$ is even and thus we now only need to consider the case of $n$
odd. Example 3 shows that $P_4$ does {\em not} divide $Q_3$.

 In the next two sections we show that for all odd $n$ with $n\geq 5, P_4$
divides $Q_n$. We first, in this section, prove the result for $n=5$. The strategy is to find a subgraph $G$ of $Q_5$, show that $G$ divides $Q_5$, and then show that $P_4$ divides $G$. In the next section we deduce the general case. 
 
 \begin{figure}[h]
 \begin{center}
 \includegraphics [scale=0.9] {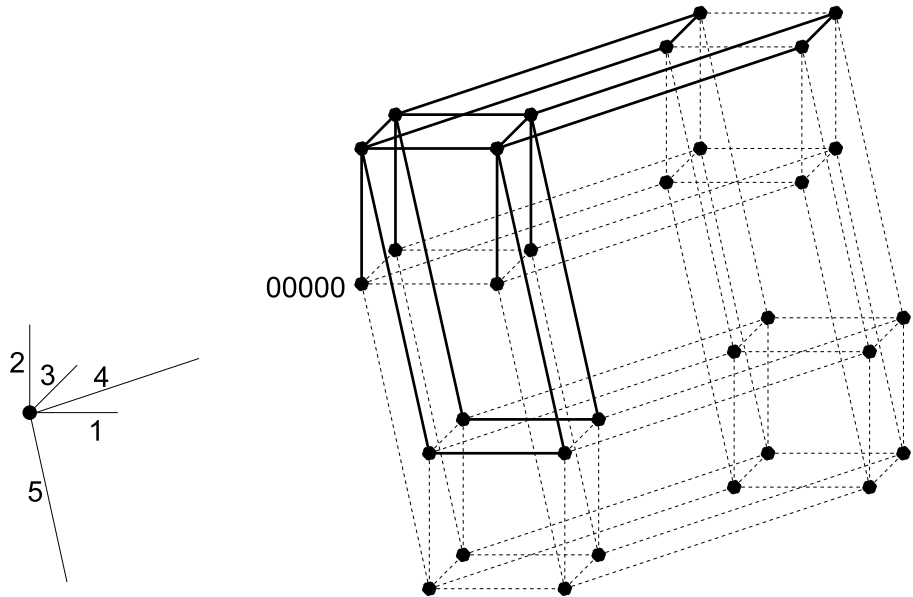}
 \caption{\label{fig:G}  $Q_5$ and the subgraph $G$}
 \end{center}
 \end{figure}
 
 We define $G$ as follows (see figure \ref{fig:G}).  First, some
 notation.  For $b, c\in \{0, 1\}, Q_5^{(***bc)}$ denotes the $3$-cube
 induced by the vertices $x_1x_2x_3x_4x_5$
 with $x_4=b$ and $x_5=c$.  If $a\in \{0,1\}$ $Q_5^{(**abc)}$ is the
 $2$-cube induced by the vertices with $x_3=a,x_4=b,$ and $x_5=c$.  We
 take $G$
 to be the union of $(1): Q_5^{(***00)}$, with the edges of $Q_5^{(*0*00)}$
 deleted; $(2)$: $Q_5^{(***10)}$ with all edges deleted except for $\langle
 01010, 01110\rangle$ and
  $\langle 11010,11110\rangle$; $(3):  Q_5^{(***01)}$ with all edges
 deleted except for $\langle 01101, 11101\rangle$ and $\langle 01001,
 11001\rangle$;
 $(4)$:  the $4$ matching edges between $Q_5^{(*1*00)}$ and
 $Q_5^{(*1*10)}$; and $(5)$  the $4$ matching edges between
 $Q_5^{(*1*00)}$  and $Q_5^{(*1*01)}$.  Thus $|\,E(G)\,|=20$.  Since
 $|\,E(Q_5)\,|=5\cdot 2^4=80$, we must exhibit $80/20=4$ copies of
 $E(G)$ that partition $E(Q_5)$.

 \begin{lemma}\label{Le:  GdivQ5}
 $G$ divides $Q_5$.  In fact, $E(G)$ is a fundamental set for $Q_5$.
 \end{lemma}
 
\begin{figure}[h]
 \begin{center}
 \includegraphics [scale=0.5] 
 {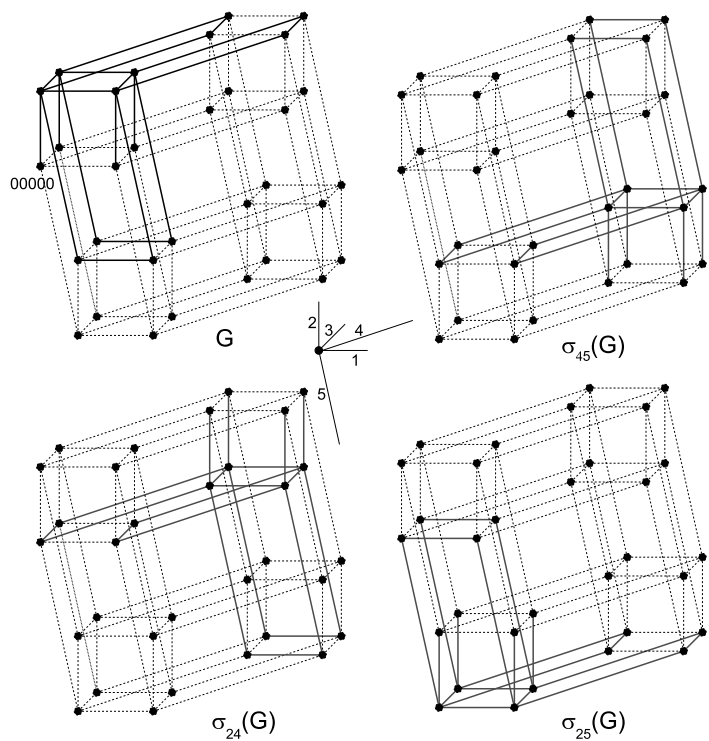}
 \caption{\label{fig:Gpart} $E(G)$ is a fundamental set for $Q_5$}
 \end{center}
 \end{figure}

 \noindent {\em Proof.} 
 By direct inspection of figure \ref{fig:Gpart} the group of translations ${\cal G}=\{id,
 \sigma_{24}, \sigma_{25}, \sigma_{45}\}$, applied to $E(G)$,
 partitions $E(Q_5)$.
%
\hfill  $\Box$
 
 \begin{figure}[h]
 \begin{center}
 \includegraphics [scale=0.9] 
 {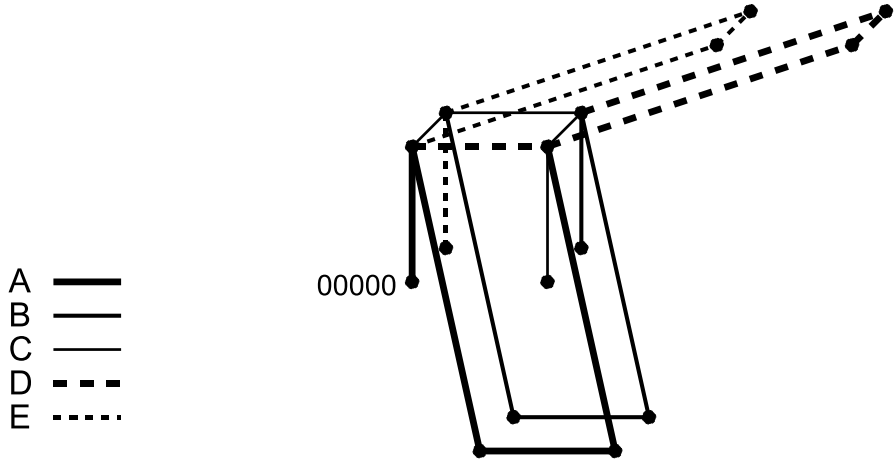}
 \caption{\label{fig:P} $P_4$ divides $G$}
 \end{center}
 \end{figure}
 
 \begin{lemma}\label{le:  P4divG}
 $P_4$ divides $G$.
 \end{lemma}
 
 \noindent  {\em Proof.} It is easiest to describe the paths by their
 starting points and direction sequences (see figure \ref{fig:P}).
 \vspace{.2in}
 
 $\begin{array}{ccc}
 \rm{Path} & \rm{Starting\,Point}  & \rm{Direction\, Sequence}\\
 A         &   00000  &  2, 5, 1, 5\\
 B         &   10100  &  2, 5, 1, 5\\
 C         &   10000  &  2, 3, 1, 3 \\
 D         &   01000  &  1, 4, 3, 4 \\
 E         &   00100  &  2, 4, 3, 4
 
  \end{array}$
 
 \hfill $\Box$
 
 \begin{corollary}\label{co:  P4divQ5}
 $P_4$ divides $Q_5$.
 \end{corollary}
 
 \noindent {\em Proof.}  This follows immediately from the previous two
 lemmas.  \hfill  $\Box$
 
 \section{$P_4$ divides $Q_n,$ for $n$ odd, $n\geq 5$}

 Let us write $Q_5$ as $Q_5=Q_3\Box Q_2=Q_3\Box C_4$.  Let
 $G_0=Q_5^{(***00)}, G_1=Q_5^{(***10)}, G_2=Q_5^{(***11)},
 G_3=Q_5^{(***01)}$.
  For $i\in \{0, 1, 2, 3\}$ let $\pi_i$ be the canonical mapping from
 $G_i$ to $Q_3$.\\
 \vspace{0.1in}
 
 \noindent $\ast$  From the decomposition of $Q_5$ by $P_4$ we have a
 coloring $c:Q_5 \longrightarrow \{1, 2, \ldots , 20\}$ of the edges of
 $Q_5$
 such that for any $i\in \{1, 2, \ldots , 20\}$ the set of edges of
 $Q_5$ colored $i$ induces a $P_4$.
 \vspace{0.1in}
 
 \noindent $\ast$  Consider now $Q_3\Box C_{4k}$ for some $k\geq 1.$
 Let $G_0^{\,\prime}, \ldots ,  G_{4k-1}^{\,\prime} \simeq Q_3$.  Let
 $\pi_{i^{\,\prime}}^{\,\prime}$
 be the canonical mapping from
 $G^{\,\prime}_{i^{\,\prime}}\longrightarrow Q_3$ for $i^{\,\prime}\in
 \{0,1, \ldots , 4k-1\}$.\\
 The edges of $Q_3\Box C_{4k}$ are\\
 \indent\indent  Case A:  the edges of $G^{\,\prime}_{i^{\,\prime}}$
 for any $i^{\,\prime} \in \{0, 1, \ldots , 4k-1\}$.\\
 \indent\indent  Case B:  for any $i^{\,\prime}\in \{0, 1, \ldots, 4k-1\}$ the
 edges $\langle x^{\,\prime}, y^{\,\prime}\rangle$ for $x^{\,\prime}\in
 G^{\,\prime}_{i^{\,\prime}},$
 $y^{\,\prime}\in G^{\,\prime}_{j^{\,\prime}}$, where
 $|\,j^{\,\prime}-i^{\,\prime}\,| \equiv 1\pmod {4k}$ and
 $\pi_{i^{\,\prime}}(x^{\,\prime})=\pi_{j^{\,\prime}}(y^{\,\prime})$.
 \vspace{0.1in}
 
 \noindent $\ast$  Let $\theta$ be the mapping from $Q_3\Box C_{4k}
 \longrightarrow Q_5$ defined by:  for any $x^{\,\prime}\in
 G^{\,\prime}_{i^{\,\prime}}, \theta (x^{\,\prime})=x$ where $x$ is the
 element of $G_i$,
 with $i\equiv i^{\,\prime}\pmod{4}$ such that
 $\pi_i(x)=\pi_{i^{\,\prime}}(x^{\,\prime})$.  (Note that $\theta$ is not a
 one-to-one mapping.)

 \begin{proposition}\label{Prop:  edges of Q5}
 If $\langle x^{\,\prime}, y^{\,\prime}\rangle$ is an edge of $Q_3\Box
 C_{4k}$ then $\langle
 \theta(x^{\,\prime}),\theta(y^{\,\prime})\rangle$ is an edge of $Q_5$.
 \end{proposition}
 
 \noindent {\em Proof.}  \\
 \indent\indent \underline{Case A}  \\
 \indent\indent\indent $\langle x^{\,\prime}, y^{\,\prime}\rangle\in
 G^{\,\prime}_{i^{\,\prime}}$ for some $i^{\,\prime}$.  Then let
 $i\equiv i^{\,\prime}\pmod {4}$.  By the definition of $\theta, \theta
 (x^{\,\prime})\in G_{i},
 \theta (y^{\,\prime})\in G_{i}$.  This implies that $\theta
 (x^{\,\prime})$ and $\theta (y^{\,\prime})$ are adjacent.
 \vspace{0.1in}
 
 \indent\indent \underline{Case B}  \\
 \indent\indent\indent  Assume $x^{\,\prime}\in
 G^{\,\prime}_{i^{\,\prime}}, y^{\,\prime} \in
 G^{\,\prime}_{j^{\,\prime}}$ with $|\,j^{\,\prime}
 -i^{\,\prime}\,|\equiv 1\pmod{4k}.$   We have
 $\pi^{\,\prime}_{i^{\,\prime}}(x^{\,\prime})=\pi^{\,\prime}_{j^{\,\prime}}
 (y^{\,\prime}).$
  Then $\theta(x^{\,\prime})\in G_i$ and $\theta(y^{\,\prime})\in G_j$
 where $|j-i|\equiv 1 \pmod{4}$ since
 $|\,j^{\,\prime}-i^{\,\prime}\,|\equiv 1\pmod{4}$ implies that
 $|\,j-i\,|\equiv 1\pmod{4}$.  Furthermore
 $$\pi_i(\theta(x^{\,\prime}))\stackrel{\rm def\, \rm of
 \theta}{=}\pi^{\,\prime}_i(x^{\,\prime})\stackrel{\rm
 edge}{=}{\pi^{\,\prime}_j(y^{\,\prime})}\stackrel{\rm def\, \rm
 of\theta}{=}{\pi_j(\theta(y^{\,\prime})}).$$
 Thus there exists an edge between $\theta(x^{\,\prime})$ and
 $\theta(y^{\,\prime})$  \hfill  $\Box$
 \vspace{0.2in}
 
 \begin{figure}[h]
 \begin{center}
 \includegraphics [scale=0.5] 
 {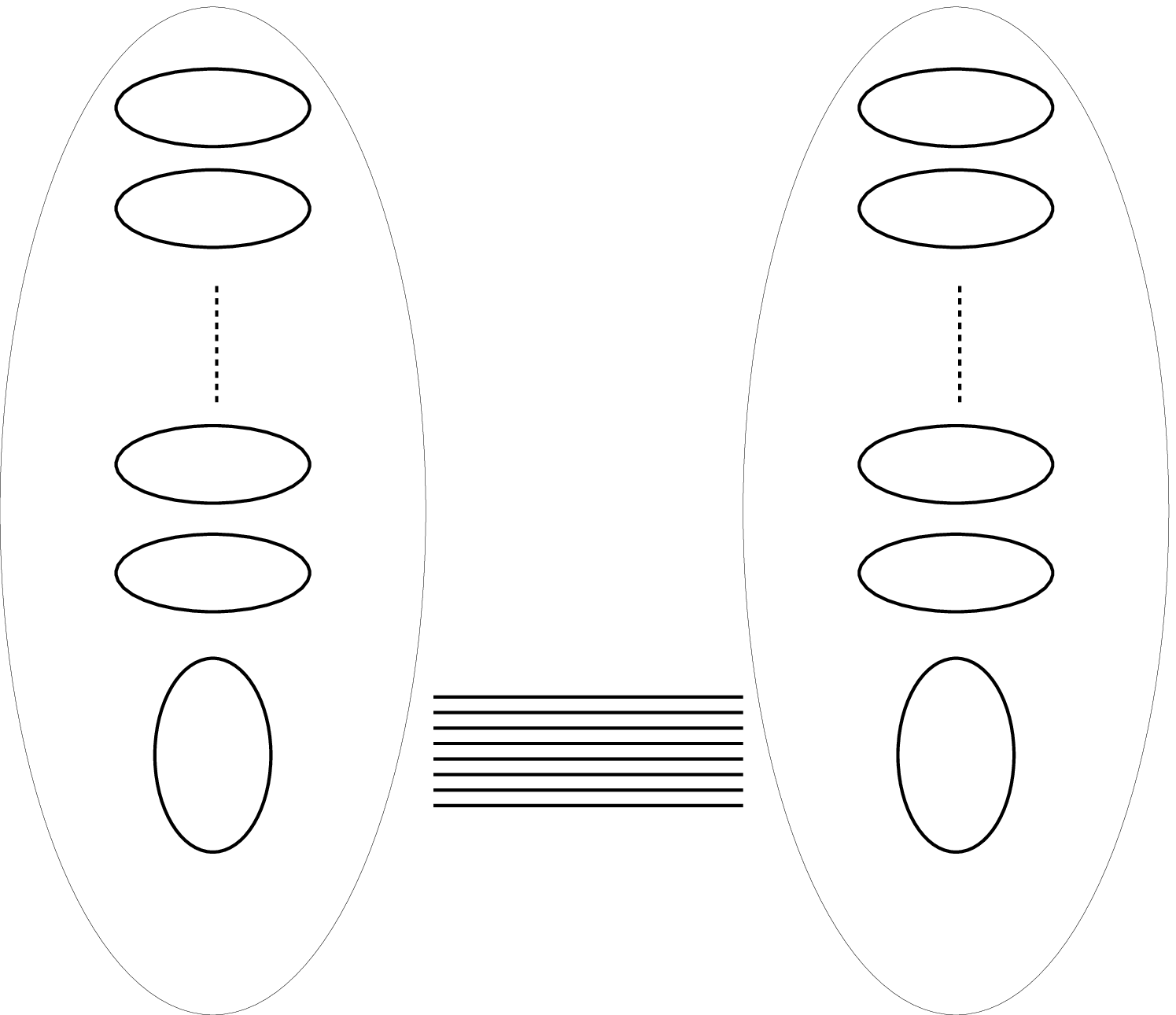}
 \caption{\label{fig:G1} Decomposition of  $Q_{2k+1}$}
 \end{center}
 \end{figure}

 \begin{definition}
 Consider the coloring $E(Q_3\Box
C_{4k})\stackrel{c^{\,\prime}}{\longrightarrow}
 \{1, 2, \ldots , 20\}$ of the edges of $Q_3\Box C_{4k}$ defined by
 $c^{\,\prime}(\langle x^{\,\prime}, y^{\,\prime}\rangle)=c(\langle \theta(x^{\,\prime}),
 \theta(y^{\,\prime})\rangle).$
 \end{definition}
 
 \begin{lemma}\label{le:C2m}
 For any $i\in \{1, 2, \ldots , 20\}$ the set of edges of $Q_3\Box C_{4k}$ such
 that $c^{\,\prime}(x^{\,\prime}, y^{\,\prime})=i$ is a set of disjoint paths of
 length $4$. Therefore $P_4$ divides $Q_3\Box C_{4m}$ for all $m\geq 1$.
 \end{lemma}
 
 \noindent {\em Proof.}
 By definition of $c^{\,\prime}$, for any vertex $x^{\,\prime}$ of
 $Q_3\Box C_{4k}$ the number of edges incident to $x^{\,\prime}$
 colored $i$ by $c^{\,\prime}$ is the number of edges incident to
 $\theta(x^{\,\prime})$ colored $i$ by $c$.  Therefore this number is
 $\leq 2$.  Furthermore, there is no cycle colored $i$ in $Q_3\Box
 C_{4k}$ because the image by $\theta$ of this cycle would be a cycle
 of $Q_5$ colored $i$ with $c$.  Therefore the set of edges colored $i$
 by $c^{\,\prime}$ is a forest and more precisely, because of the
 degree, a set of disjoint paths.
 
 \noindent Notice that the image by $\theta$ of a path colored $i$ is a
 path of $Q_5$ of the same length (because of the degree of the
 endpoints of the paths).  Therefore all the paths are of length $4$.
 \hfill  $\Box$

 \begin{figure}[h]
 \begin{center}
 \includegraphics [scale=0.5] 
{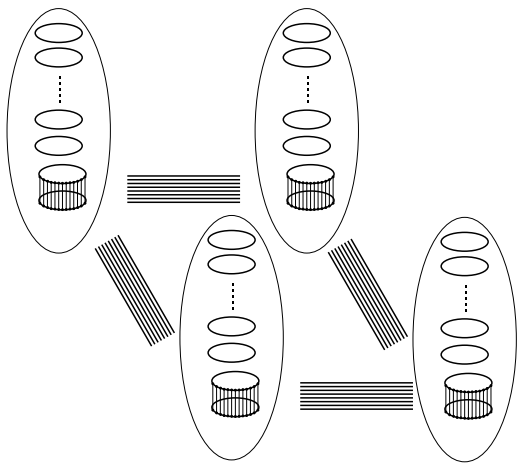}
 \caption{\label{fig:dec} Decomposition of  $Q_{2k+3}$}
 \end{center}
 \end{figure}

 \begin{theorem}\label{P4divQn}
  For $n\geq 4$, $P_4$ divides $Q_n$.
 \end{theorem}
 
 \noindent {\em Proof.}
 \noindent
 If $n$ is even, the result is true by Corollary \ref{P2kdivQn}. If $n=5$ then we are done by Corollary \ref{co:  P4divQ5}.
 Consider $Q_{2k+3}$, for $k\geq 2$.  $Q_{2k+3}=Q_{2k+1}\Box Q_2$.  $E(Q_{2k})$
 can be decomposed into $k$ cycles of length
 $2^{2k}$ (Hamiltonian cycles) by Aubert and Schneider \cite{AS}. Let $D$ be one
 of these cycles. The edges of $Q_{2k+1}$ are the edges of the two copies of
 $Q_{2k}$ and a matching. But every vertex of $Q_{2k}$ appears exactly once in
 $D$ so $E(Q_{2k+1})$ can be decomposed into $2(k-1)$ cycles of length $2^{2k}$
 and $D\Box Q_1\simeq C_{2^{2k}}\Box Q_1$ (see figure \ref{fig:G1}).

 Every vertex of $Q_{2k+1}$ appears once in $D\Box Q_1$,
 thus, for the same reason, $E(Q_{2k+3})$ can be decomposed into $8(k-1)$ cycles
 of length $2^{2k}$ and $D\Box Q_1\Box Q_2 \simeq  C_{2^{2k}}\Box Q_1\Box Q_2
 \simeq C_{2^{2k}}\Box Q_3$ (see figure \ref{fig:dec}).\\
 
 Since $k\geq 2$, $\frac{2^{2k}}{4}$ is an integer strictly greater than $1$  so
 the cycles of length $2^{2k}$ are divisible by $P_4$.  By Lemma \ref{le:C2m},
 $P_4$ divides $C_{2^{2k}}\Box Q_3$, and
 $P_4$ divides $E(Q_{n})$ for any odd $n\geq5$.
  \hfill   $\Box$

 \section{$Q_{2^k}$ has a fundamental Hamiltonian cycle.}

   We shall
 describe walks in the hypercube by specifying the starting vertex (generally
 $\emptyset$) and the sequence of edge directions.
 
 \indent\indent It is well-known that the $n$-dimensional hypercube
 $Q_n$ is Hamiltonian, and in fact has many Hamiltonian cycles.  Aubert
 and Schneider \cite{AS}
 proved that for $n$ even, $Q_n$ has an edge decomposition into Hamiltonian
 cycles.  However, their construction is technical.  In contrast, in this last
 section we shall
 prove that for $n=2^k$, there is a single Hamiltonian cycle $C$
  such that $E(C)$ is a fundamental set for $Q_n$.
 
 By $G_1\Box G_2$ we denote the Cartesian product of the graphs $G_1$ and $G_2$.
 We will start with two easy results about Cartesian product of graphs.
 \begin{lemma} \label{le:DivCartsquare}
  Assume that $\{C^1,C^2,\dots, C^p  \}$ is  an edge decomposition in Hamiltonian cycles of a graph $G$.  Then $\{C^1 \Box C^1 ,C^2 \Box C^2,\dots,C^p \Box C^p\}$ is an edge decomposition of $G \Box G$.
 \end{lemma}
 \noindent  {\em Proof.}
 Let $(x_1, x_2)$ and $(y_1, y_2)$ be adjacent
 in $G \Box G $.  Then either  $x_1$ and $y_1$ are adjacent in
 $G$ and $x_2=y_2$ or $x_1=y_1$ and $x_2$ and $y_2$ are adjacent
 in $G$. By symmetry, it is sufficient to consider the first case. Let $i$ be such that $\langle x_1,y_1 \rangle\in E(C^i)$.  Then since $C^i$ is Hamiltonian $x_2=y_2 \in V(C^i)$;  thus  $\langle(x_1, x_2),(y_1, y_2)\rangle\in E(C^i \Box C^i)$. Conversely  $\langle(x_1, x_2),(y_1, y_2)\rangle\in E(C^j  \Box C^j)$ implies $\langle x_1,y_1 \rangle\in E(C^j)$ since $x_2=y_2$; thus $j=i$. Therefore the $C^j \Box C^j$ 's are disjoint and the conclusion follows.
  \hfill   $\Box$
  
  \begin{lemma}  \label{AutProd}
 Let $G_1$ and $G_2$ be any two graphs, and for $i=1, 2$ let $\phi_i\in
 \, {\rm Aut}\,(G_i)$.  Define $(\phi_1, \phi_2): G_1\Box G_2
 \longrightarrow G_1\Box G_2$ by $(\phi_1, \phi_2)((x, y))=(\phi_1(x),
 \phi_2(y))$.  Then
 $(\phi_1, \phi_2)\in\, {\rm Aut}\,(G_1\Box G_2).$
 \end{lemma}
 
 \noindent {\em Proof.}  Let $(x_1, x_2)$ and $(y_1, y_2)$ be adjacent
 in $G_1\Box G_2$.  Then either (1)  $x_1$ and $y_1$ are adjacent in
 $G_1$ and $x_2=y_2$ or (2)  $x_1=y_1$ and $x_2$ and $y_2$ are adjacent
 in $G_2$.  We must show that $(\phi_1, \phi_2)(x_1, x_2)$ and
 $(\phi_1, \phi_2)(y_1, y_2)$ are adjacent in $G_1\Box G_2$.  By
 symmetry, it is sufficient to prove this for case (1).  But then since
 $\phi_1 \in \,{\rm Aut}\,(G_1)$,
 $\phi_1(x_1)$ and $\phi_1(y_1)$ are adjacent in $G_1$, and since $x_2=y_2,\,
 \phi_2(x_2)=\phi_2(y_2)$.  Therefore $(\phi_1, \phi_2)(x_1, x_2)$ and $(\phi_1,
 \phi_2)(y_1, y_2)$ are adjacent in $G_1\Box G_2$.  Conversely if $(\phi_1, \phi_2)(x_1, x_2)=(\phi_1(x_1), \phi_2(x_2))$ and $(\phi_1,
 \phi_2)(y_1, y_2)=(\phi_1(y_1), \phi_2(y_2))$ are adjacent in $G_1\Box G_2$ then $\phi_1(x_1)=\phi_1(y_1)$ or $\phi_2(x_2)=\phi_2(y_2)$. We can assume the first case by symmetry then  $x_1=y_1$  and $x_2$ is adjacent to $y_2$ in $G_2$.
 Thus $(x_1, x_2)$ and $(y_1, y_2)$ are adjacent
 in $G_1\Box G_2$ and $(\phi_1, \phi_2)\in
 \,{\rm Aut}\,(G_1\Box G_2)$.
  \hfill   $\Box$
\\
 The starting point of the theorem of Aubert and Schneider is an earlier result of G. Ringel \cite {Ring} who proved that for  $n=2^k$, $Q_n$ has an edge decomposition into Hamiltonian
 cycles. His proof is by induction on $k$. Let us recall the induction step.
Let $m=2^n$. Let $\theta$ be the mapping from $\{1,\dots,n\}$ to $\{n+1,\dots,2n\}$ defined by $\theta(i)=i+n$.
 Let $C$ be a Hamiltonian cycle of $Q_n$ then we can construct $\Phi(C)$ and $\Gamma(C)$ two disjoint Hamiltonian cycles of $Q_{2n}=Q_n\Box Q_n$ such that $E(C \Box C)=E(\Phi(C))\cup E(\Gamma(C))$. Indeed fix an arbitrary vertex (say $0$) and represent $C$ by the sequence of directions $C=(c_1,\dots,c_m)$ then consider\\
$\begin{array}{llll}
 \Phi(C)=(         &c_1,...&...,c_{m-1},c_{\theta(c_1)},&\\
        &c_m,c_1,...&...,c_{m-2},c_{\theta(c_2)},&\\
          &c_{m-1},c_m,c_1,...&...,c_{m-3},c_{\theta(c_3)},&\\
        &.....&.....&\\
           &c_2,...&...,c_{m},\;c_{\theta(c_m)},&)\\
 
  \end{array}$\\
  and\\
  $\begin{array}{llll}
\Gamma(C)=(         &c_{\theta(1)},...&...,c_{\theta(m-1)},c_1,&\\
        &c_{\theta(m)},c_{\theta(1)},...&...,c_{\theta(m-2)},c_2,&\\
          &c_{\theta(m-1)},c_{\theta(m)},c_{\theta(1)},...&...,c_{\theta(m-3)},c_3,&\\
          &.....&.....&\\
           &c_{\theta(2)},...&...,c_{\theta(m)},\;c_m,&)\\
 
  \end{array}$.\\
  \begin{figure}[h]
 \begin{center}
 \includegraphics [scale=0.5] {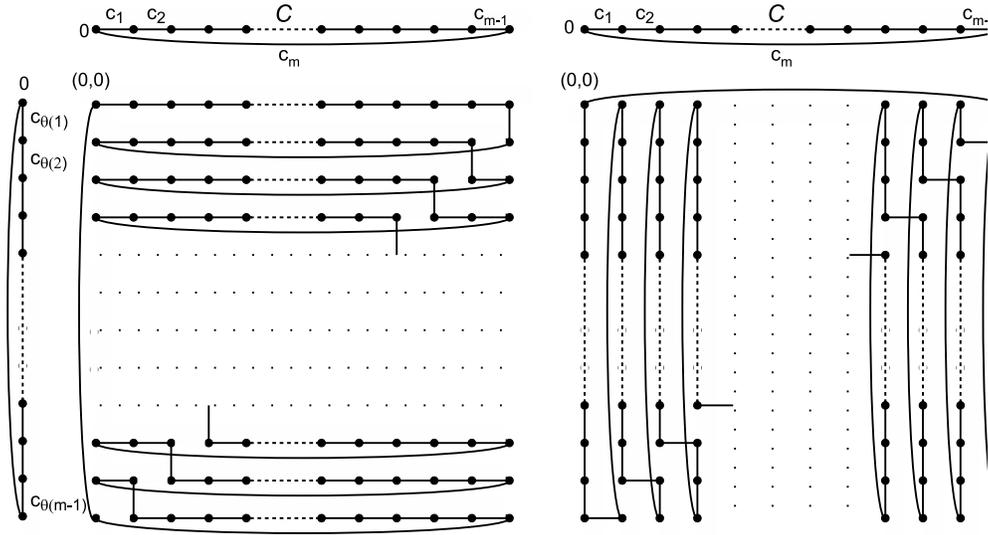}
 \caption{\label{fig:Ringel} Construction of $\Phi(C)$ and $\Gamma(C)$ from $C$}
 \end{center}
 \end{figure}
 
  It is immediate to check (see figure \ref{fig:Ringel}) that $\Phi(C)$ and $\Gamma(C)$ are disjoint and define a partition of the edges of $C \Box C$. For $n$ even let $p=n/2$ and  assume that $\{C^1,C^2,\dots,C^p\}$ is an edge decomposition of $Q_n$ in Hamiltonian cycles then as a consequence of Lemma \ref{le:DivCartsquare}, $\{\Phi(C^1),\Phi(C^2),\dots,\Phi(C^p)\}\cup \{\Gamma(C^1),\Gamma(C^2),\dots,\Gamma(C^p)\}$ is an edge decomposition of $Q_{2n}$ in Hamiltonian cycles. 
 
  \begin{theorem}  \label{theor:  FundHam}
 For any $k\geq 1, Q_{2^k}$ has a Hamiltonian cycle that is a fundamental set.
 \end{theorem}
 
\noindent  {\em Proof.}  This is trivial for $k=1$ since $Q_2=C_4$.  The desired result follows by
 induction from Ringel's construction.
  Indeed  let $n=2^k, k\geq1$  and assume that there exists  an edge  decomposition $\{C^1,C^2,\dots,C^p\}$ of $Q_n$ obtained as the translate of an Hamiltonian cycle $C^1$ under some subgroup  ${\cal E}$ of {\rm Aut} ($Q_n$). For any automorphism $\phi\in$ {\rm Aut} ($Q_n$), $(\phi,\phi)\in ${\rm Aut} ($Q_{2n}$) by Lemma {\ref{AutProd}. Furthermore if $\phi(C^1)=C^i$ then $(\phi,\phi)(\Phi(C^1))=\Phi(C^i)$ and $(\phi,\phi)(\Gamma(C^1))=\Gamma(C^i)$.
  If we consider now the permutation $\theta$ on $\{1,\dots,2n\}$ defined by $\theta(i)=i+n\: mod\:2n$ then $\rho_{\theta}(\Phi(C^i))=\Gamma(C^i)$. The conclusion follows since the subgroup of {\rm Aut} ($Q_{2n}$), isomorphic to ${\cal E}\times S_2$, defined by  ${\cal H}=\{(\phi,\phi);\phi \in \cal E\}\cup\{\rho_{\theta}\circ (\phi,\phi);\phi \in \cal E\}$ is such that $\{\Phi(C^1),\Phi(C^2),\dots,\Phi(C^p)\}\cup \{\Gamma(C^1),\Gamma(C^2),\dots,\Gamma(C^p)\}$ are the translates of  $\Phi(C^1)$ under ${\cal H}$.
  \hfill   $\Box$ 
  \begin{corollary} \label{n,m powers of 2}
  For $n$ and $m$ each a power of $2$, with $m\leq n$,  there is an $m$-cycle
 that divides $Q_n$.
 \end{corollary}
 
 \noindent {\em Proof.}  Let $m=2^p$.  By Theorem \ref{theor:  FundHam}
 $Q_m$ has a
 fundamental
 $2^p$-cycle, which therefore divides $Q_m=Q_{2^p}$.  Since $m$ and $n$ are each
 powers of two, $m$ divides $n$. Hence
  by Proposition \ref{QkdivQn} and Lemma \ref{Lemmadiv}, this cycle divides
 $Q_n$.  \hfill   $\Box$
%
%
%
%
%
 
 \section{Summary of Results}
 
 1.   For $k$ odd, if $P_k$ is a path on $k$ edges that divides $Q_n$, then $k$ divides $n$.
 (Proposition \ref{kdivn}) \\

 \noindent 2.  If $k$ divides $n$, {\em any} tree on $k$ edges divides $Q_n$.  (Corollary \ref{TdivQn})  \\
 3.  If $k$ divides $n$ and $k<n$ then $P_{2k}$ divides $Q_n$.  (Corollary
 \ref{P2kdivQn})\\
 4.  If $n$ is even and $j<n$ then $P_{2^j}$ divides $Q_n$.
 (Proposition \ref{Prop:  P_2^j})\\
 5.  For $k=2n$ there is a $k$-cycle which is a fundamental set for $Q_n$ when
 $n$ is even. (Proposition \ref{cycle}) \\
 6.  For $n=$ a power of $2$, there is a Hamiltonian cycle which is a
fundamental
 set for $Q_n$. (Theorem \ref{theor:  FundHam})\\
 7.  For $n=$ a power of $2$ and $m=$ a power of $2$, with $m\leq n$, there is
an
 $m$-cycle that divides $Q_n$.  (Corollary \ref{n,m powers of
 2})\\
 8.  For $n\geq 4$, $P_4$ divides $Q_n$. (Theorem \ref{P4divQn})\\
 9.   $Q_k$ is a fundamental set for $Q_n$ if and only if $k$ divides $n$.
 (Proposition \ref{QkdivQniff kdivn})\\
 10.  For $k\geq 3, P_{2^j}$ does {\em not} divide $Q_{2k+1}.$ 
 (Proposition \ref{doesNotdiv})

\end{document}